\documentclass[12pt]{article}
\usepackage{amsmath}
\usepackage{graphicx}
\usepackage{CJK}
\usepackage{multirow}
\usepackage{makecell}

\usepackage{subfigure}
\usepackage[justification=centering]{caption}
\input{amssym.tex}

\def\bc{\begin{center}}       \def\ec{\end{center}}
\def\ba{\begin{array}}        \def\ea{\end{array}}
\def\be{\begin{equation}}     \def\ee{\end{equation}}
\def\bea{\begin{eqnarray}}    \def\eea{\end{eqnarray}}
\def\beaa{\begin{eqnarray*}}  \def\eeaa{\end{eqnarray*}}

\def\mathbb{\Bbb}

\begin{document}
\baselineskip 18pt
\centerline{\Large\bf Limit cycles bifurcating from the quasi-homogeneous polynomial    }
\vskip 0.2 true cm
\centerline{\Large\bf centers  of weight-degree 2 under non-smooth perturbations }

\vskip 0.3 true cm

\centerline{\bf  Shiyou Sui$^{\dag}$\footnote{Author for correspondence. E-mail: sui\_shiyou@163.com (S. Sui), changyongkang@buaa.edu.cn (Y. Zhang), libaoyi1123@euou.com(B. Li).}, Yongkang Zhang$^{\ddag}$, Baoyi Li$^{\ddag}$}
\centerline{$^\dag$ School of Science, Tianjin University of Commerce,}
\centerline{Tianjin 300134, The People's Republic of China}
\centerline{$^\ddag$ School of Mathematical Sciences, Tianjin Normal University,}
\centerline{Tianjin 300387, The People's Republic of China}

\vskip 0.2 true cm
\noindent{\bf Abstract} We investigate the maximum number of limit cycles bifurcating from the period annulus of a family of cubic polynomial differential centers when it is perturbed inside the class of all  cubic piecewise smooth polynomials. The family considered is the unique family of weight-homogeneous polynomial differential systems of weight-degree 2 with a center. When the switching line is  $x=0$ or $y=0$, we obtain the sharp bounds of the number of limit cycles for the perturbed systems by using the first order averaging method.
Our results indicate that non-smooth systems can have more limit cycles than smooth ones, and the switching lines play an important role in the dynamics of non-smooth systems.


\vskip 0.2 true cm
\centerline {\bf \large 1 Introduction and statement of the main results}
\vskip 0.2 true cm

A periodic orbit of a differential system which is isolated in the set of all periodic orbits of the system is a {\it limit cycle}. One of the main goals in the qualitative theory of differential system  is the study of limit cycles of planar polynomial system
$$
\begin{cases}
\displaystyle\frac{{\rm d}x}{{\rm d}t}=P(x,y),\\
\displaystyle\frac{{\rm d}y}{{\rm d}t}=Q(x,y),
\end{cases}\eqno{(1.1)}
$$
where $P$ and $Q$ are polynomials with real coefficients, which is known as Hilbert's 16th problem \cite{D00,S98}. The {\it degree} of system $(1)$ is the maximum of the degrees of the polynomials  $P$ and $Q$. A singular point $p$ of system $(1)$ is a {\it center} if there is a neighborhood of $p$ fulfilled by periodic orbits. The { \it period annulus} of a center  is the region fulfilled by all the periodic orbits  surrounding the centers. We say a center at the origin of coordinates is {\it global} if its period annulus is $\mathbb{R}^2\setminus\{(0,0)\}$. Many authors  have studied the number of limit cycles bifurcating from the period annuluses, which is related to weaken Hilbert's 16th problem \cite{L06,L03}.

The polynomial differential  system $(1.1)$ is {\it quasi-homogeneous} if there exist $s_1,s_2,d\in\mathbb{N}$ such that for arbitrary $\lambda\in\mathbb{R}^+=\{\lambda\in\mathbb{R}:\lambda>0\}$,
$$P(\lambda^{s_1}x,\lambda^{s_2}y)=\lambda^{s_1-1+d}P(x,y),~~~Q(\lambda^{s_1}x,\lambda^{s_2}y)=\lambda^{s_2-1+d}Q(x,y).\eqno{(1.2)}$$
We call $s=(s_1,s_2)$ the {\it weight exponent} of system $(1.1)$ and $d$ the {\it weight degree} with respect to the weight exponent $s$. Particularly, if $s_1=s_2=1$, then system $(1.1)$ is the classical homogeneous polynomial differential system of degree $d$. The quasi-homogeneous polynomial differential systems have been intensively studied by a great deal of  authors from  different points of view, such as normal forms \cite{GLP13}, integrability \cite{G03}, center \cite{GLV17},  limit cycles \cite{LLYZ09}.

The classification of all centers of planar weight-homogeneous polynomial differential systems up to weight-degree 4 can be found in \cite{LP09}. They proved that the unique family of quasi-homogeneous polynomial differential system having center with weight-degree 2 is
$$
\begin{cases}
\displaystyle\frac{{\rm d}x}{{\rm d}t}=ax^2+by,\\
\displaystyle\frac{{\rm d}y}{{\rm d}t}=cx^3+dxy,
\end{cases}\eqno{(1.3)}
$$
where
$$(d-2a)^2+8bc<0.$$
The weight-exponent of system $(1.3)$ is $(s_1,s_2)=(1,2)$. From \cite{GGL15,LLYZ09}, we know that system $(1.3)$ has a global center at the origin of coordinates. In \cite{LL16}, the authors studied the number of limit cycles bifurcating from the center of system $(1.3)$ under perturbation of smooth polynomials of degree 3, and they solved the weak Hilbert's 16th problem for this case.

Stimulated by discontinuous phenomena in the real world, there are lots of works in mechanics, electrical engineering and the theory of automatic control which are described by non-smooth systems \cite{dBCK08}.  A big interest has appeared for studying the limit cycles of discontinuous differential systems \cite{LM14,LCZ17,LH10,YZ18}. Our objective is to study the   maximum number of limit cycles which can bifurcate from the periodic orbits of the quasi-homogeneous centers $(1.3)$ perturbed inside  discontinuous cubic polynomials.  More precisely we consider the following two non-smooth differential systems
$$
\begin{cases}
\displaystyle\frac{{\rm d}x}{{\rm d}t}=ax^2+by+\varepsilon p_1(x,y),\\
\displaystyle\frac{{\rm d}y}{{\rm d}t}=cx^3+dxy+\varepsilon q_1(x,y),
\end{cases}\eqno{(1.4)}
$$
with
$$
p_1(x,y)=\begin{cases}\sum\limits _{i+j=0}^3a_{i,j}^+x^iy^j,~~x\geq0,\\
\sum\limits _{i+j=0}^3a_{i,j}^-x^iy^j,~~x<0,
\end{cases}
~~~
q_1(x,y)=\begin{cases}\sum\limits _{i+j=0}^3b_{i,j}^+x^iy^j,~~x\geq0,\\
\sum\limits _{i+j=0}^3b_{i,j}^-x^iy^j,~~x<0,
\end{cases}$$
and
$$
\begin{cases}
\displaystyle\frac{{\rm d}x}{{\rm d}t}=ax^2+by+\varepsilon p_2(x,y),\\
\displaystyle\frac{{\rm d}y}{{\rm d}t}=cx^3+dxy+\varepsilon q_2(x,y),
\end{cases}\eqno{(1.5)}
$$
with
$$
p_2(x,y)=\begin{cases}\sum\limits _{i+j=0}^3c_{i,j}^+x^iy^j,~~y\geq0,\\
\sum\limits _{i+j=0}^3c_{i,j}^-x^iy^j,~~y<0,
\end{cases}
~~~
q_2(x,y)=\begin{cases}\sum\limits _{i+j=0}^3d_{i,j}^+x^iy^j,~~y\geq0,\\
\sum\limits _{i+j=0}^3d_{i,j}^-x^iy^j,~~y<0,
\end{cases}$$
where $a_{i,j}^{\pm},~b_{i,j}^{\pm},~c_{i,j}^{\pm},~d_{i,j}^{\pm}$ are any real constants, and $\varepsilon$ is a small parameter.

Applying the  averaging  theory derived in \cite{WZ18}, we bound the number of limit cycles bifurcating from the periodic orbits of system $(1.2)$. The main results are follows.

\noindent{\bf Theorem 1.1}. Consider the perturbed system $(1.4)$, the following two statement hold.\\
(a) There are at most 7 limit cycles for system $(1.4)$ by using the first order averaging function.\\
(b) There exist perturbations such that system $(1.4)$ has exactly 7 limit cycles.

\noindent{\bf Theorem 1.2}. For the perturbed system $(1.5)$,  we have the  following two statement.\\
(a) There are at most 3 limit cycles for system $(1.5)$ by using the first order averaging function.\\
(b) There exist perturbations such that system $(1.5)$ has exactly 3 limit cycles.

In some sense we extend the work done by Llibre and  de Moraes \cite{LL16} for the continuous polynomial perturbations to the discontinuous ones with the straight line of discontinuity $x=0$ or $y=0$. Recall that the perturbations of the periodic orbits of the quasi-homogeneous centers $(1.3)$ inside the class of continuous cubic polynomial differential systems produce at most 3 limit cycles. Comparing the results obtained in Theorem  1.1 with the continuous case, this work shows that the discontinuous systems can have 4 more limit cycles surrounding the origin than the continuous systems when we perturb the quasi-homogeneous centers $(1.3)$. Moreover, comparing the results of Theorem 1.1 with Theorem 1.2, the switching lines play an important role in the dynamics of non-smooth differential systems.

The rest of this paper is organized as follows. In Section 2, we give some preliminary results. Theorem 1.1 and 1.2 are proved in Section 3 and Section 4, respectively.

\vskip 0.2 true cm
\centerline {\bf \large 2 Preliminary results}
\vskip 0.2 true cm

This section is devoted to present some  preliminary tools needed to prove our main results. Firstly, we introduce the averaging theory for discontinuous differential systems proved by \cite{WZ18}.

Consider a piecewise  differential equation
$$\frac{{\rm d}r}{{\rm d}\theta}=F_0(\theta,r)+\varepsilon F_1(\theta,r)+O(\varepsilon^2),\eqno{(2.1)}$$
with
$$F_i(\theta,r)=\begin{cases}F_i^+(\theta,r), ~~{\text{if}} ~~\alpha<\theta<\alpha+\pi,\\
F_i^-(\theta,r), ~~{\text {if}} ~~\alpha-\pi<\theta<\alpha,
\end{cases} $$
where $F_i^{\pm}:[\alpha-\pi,\alpha+\pi]\times(0,\rho^*)\rightarrow \mathbb{R}$ are analytical functions $2\pi-$periodic in the variable $\theta$ for $i=1,2$, and $\varepsilon\in(-\varepsilon_0,\varepsilon_0)$ with $\varepsilon_0$ a small positive real number.

Denote by $r(\theta;z,\varepsilon)$ the solution of system $(2.1)$ with the initial condition $r(\alpha)=z$. And we use the notation $r(\theta;\alpha+\pi,z,\varepsilon)$ to denote the solution of system $(2.1)$ satisfying the initial condition $r(\alpha+\pi)=z$.

Suppose that the unperturbed equation of $(2.1)$ i.e.
$$\frac{{\rm d}r}{{\rm d}\theta}=F_0(\theta,r),\eqno{(2.2)}$$
has a family of periodic orbits of period $2\pi$, which are filled with a region of $[\alpha-\pi,\alpha+\pi]\times(0,\rho^*)$. Let $r_0(\theta;z)$ be a solution of the unperturbed system $(2.2)$ satisfying $r_0(\alpha)=z$. The solution $r_0(\theta;z)$ of $(2.2)$ can be seen as a composition of the solution $r_0^+(\theta;z)$ of the initial value problem
$$\frac{{\rm d}r}{{\rm d}\theta}=F_0^+(\theta,r),~~~~r(\alpha)=z,\eqno{(2.3)}$$
when $\theta\in[\alpha,\alpha+\pi]$, and of the solution $r_0^-(\theta;\alpha+\pi,w)$ of the initial value problem
$$\frac{{\rm d}r}{{\rm d}\theta}=F_0^-(\theta,r),~~~r(\alpha+\pi)=w:=r_0^+(\alpha+\pi,z),\eqno{(2.4)}$$
when $\theta\in[\alpha-\pi,\alpha]$. That is
$$r_0(\theta;z)=\begin{cases}r_0^+(\theta;z), ~~~~~~~~~~~~~~~~~~~~~~~~~~~~~~~~~~~~~~~~~~~~~~~~~~~\theta\in[\alpha,\alpha+\pi],\\
r_0^-(\theta-\alpha+\pi;r_0^+(\alpha+\pi,z))=r_0^-(\theta;\alpha+\pi,w), ~~\theta\in[\alpha-\pi,\alpha].
\end{cases} $$
In the following, we will use $r_0^-(\theta,z)$ to represent $r_0^-(\theta;\alpha+\pi,r_0^+(\alpha+\pi,z))$.

We define the first order averaging function $h_1:(0,\rho^*)\rightarrow\mathbb{R}$ as
$$h_1(z)=\int_{\alpha}^{\alpha+\pi}\frac{F_1^+(s,r_0^+(s;z))}{\frac{\partial{r_0^+(s;z)}}{\partial{z}}}{\rm d}s+\int_{\alpha-\pi}^{\alpha}\frac{F_1^-(s,r_0^-(s,z))}{\frac{\partial{r_0^-(s;z)}}{\partial{z}}}{\rm d}s.\eqno{(2.5)}$$
Let $r(\theta;z,\varepsilon)$ be a solution of $(2.1)$ satisfying $r(\alpha)=z$.

\noindent{\bf Theorem 2.1} (\cite{WZ18}).  Suppose that the solution $r_0(\theta;z)$ of $(2.2)$ satisfying the initial condition $r_0(\alpha)=z$ is of $2\pi-$periodic for $z\in(0,\rho^*)$. If $h_1(z)$ is not identically zero, then for each simple root $z^*$ of $h_1(z)=0$ and $\left|\varepsilon\right|\neq0$ sufficiently small, there exists a $2\pi-$periodic solution $r(\theta;\phi(\varepsilon),\varepsilon)$ of $(2.1)$ such that $r(\alpha;\phi(\varepsilon),\varepsilon)\rightarrow z^*$ as $\varepsilon\rightarrow 0$, where $\phi$ is an analytic function which satisfies $\phi(0)=z^*$.

To estimate the number of zeros of the first order averaging function, we need the following generalized Descartes Theorem proved in \cite{BZ64}.

\noindent{\bf Theorem 2.2}. Consider the real polynomial $q(x)=a_{i_1}x^{i_1}+a_{i_2}x^{i_2}+\cdots+a_{i_r}x^{i_r}$ with $0\leq i_1<i_2<\cdots<i_r$. If $a_{i_j}a_{i_{j+1}}<0$, we say that we have a variation of sign. If the number of variations of signs is $m$, then the polynomial $q(x)$ has at most $m$ positive real roots. Furthermore, always we can choose the coefficients of polynomial $q(x)$ in such a way that $q(x)$ has exactly $r-1$ positive real roots.

In order to simplify the calculation process in the following discussion, we denote by
$$M(\cos\theta,\sin\theta)=a\cos^2\theta+b\sin\theta,\eqno{(2.6)}$$
$$N(\cos\theta,\sin\theta)=c\cos^3\theta+d\cos\theta\sin\theta,\eqno{(2.7)}$$
$$f(\theta)=\cos\theta M(\cos\theta,\sin\theta)+\sin\theta N(\cos\theta,\sin\theta),\eqno{(2.8)}$$
and $$g(\theta)=\cos\theta N(\cos\theta,\sin\theta)-2\sin\theta M(\cos\theta,\sin\theta).\eqno{(2.9)}$$

\vskip 0.2 true cm
\centerline {\bf \large 3 Proof of Theorem 1.1}
\vskip 0.2 true cm

Considering the perturbed system $(1.4)$ and the weighted blow-up $x=r\cos\theta,~~y=r^2\sin\theta$, we have that
$$\frac{{\rm d}r}{{\rm d}\theta}=F_0(\theta,r)+\varepsilon F_1(\theta,r)+O(\varepsilon),\eqno{(3.1)}$$
in the standard form for applying the averaging theory of  first order described in Section 2, where
$$F_0(\theta,r)=\frac{f(\theta)}{g(\theta)}r,\eqno{(3.2)}$$
$${\begin{split}
F_1(\theta,r)&=\frac{1+\sin^2\theta}{r^3g(\theta)^2}\left[r^2N(\cos\theta,\sin\theta)p_1(r\cos\theta,r^2\sin\theta)\right.\\
&-\left.rM(\cos\theta,\sin\theta)q_1(r\cos\theta,r^2\sin\theta)\right],
\end{split}}\eqno{(3.3)}$$
and
$$p_1(r\cos\theta,r^2\sin\theta)=\begin{cases}\sum\limits _{i+j=0}^3a_{i,j}^+r^{i+2j}\cos^i\theta \sin^j\theta,~~\cos\theta\geq0,\\
\sum\limits _{i+j=0}^3a_{i,j}^-r^{i+2j}\cos^i\theta \sin^j\theta,~~\cos\theta<0,
\end{cases}$$
$$q_1(r\cos\theta,r^2\sin\theta)=\begin{cases}\sum\limits _{i+j=0}^3b_{i,j}^+r^{i+2j}\cos^i\theta \sin^j\theta,~~\cos\theta\geq0,\\
\sum\limits _{i+j=0}^3b_{i,j}^-r^{i+2j}\cos^i\theta \sin^j\theta,~~\cos\theta<0.
\end{cases}$$

Hence, the equation $(3.1)$ with $\varepsilon=0$ can be written into the form
$$\frac{{\rm d}r}{{\rm d}\theta}=\frac{f(\theta)}{g(\theta)}r,\eqno{(3.4)}$$
The solution of equation $(3.4)$  with initial value $r(-\frac{\pi}{2})=z$ is
$$r_0(\theta;z)=u(\theta)z,\eqno{(3.5)}$$
where $u(\theta)={\rm exp}\left(\int_{-\frac{\pi}{2}}^\theta\frac{f(\tau)}{g(\tau)}{\rm d}\tau\right)$. It is easy to check that $F_1(\theta,r)$ is analytical function $2\pi-$periodic in the variable $\theta$. Then, by $(2.4)$, we know that the averaging function for system $(3.1)$ is
$${\begin{split}
h_1(z)&=\int_{-\frac{\pi}{2}}^{\frac{3\pi}{2}}\frac{F_1(s,u(s)z)}{u(s)}{\rm d}s\\
&=\int_{-\frac{\pi}{2}}^{\frac{\pi}{2}}\frac{1+\sin^2\theta}{u(s)^4z^3g(s)^2}\left[u(s)^2z^2N(\cos s,\sin s)\sum\limits_{i+j=0}^3a_{i,j}^+u(s)^{i+2j}z^{i+2j}\cos^is\sin^js\right.\\
&~~~\left.-u(s)zM(\cos s,\sin s)\sum\limits_{i+j=0}^3b_{i,j}^+u(s)^{i+2j}z^{i+2j}\cos^is\sin^js\right]{\rm d}s\\
&+\int_{\frac{\pi}{2}}^{\frac{3\pi}{2}}\frac{1+\sin^2\theta}{u(s)^4z^3g(s)^2}\left[u(s)^2z^2N(\cos s,\sin s)\sum\limits_{i+j=0}^3a_{i,j}^-u(s)^{i+2j}z^{i+2j}\cos^is\sin^js\right.\\
&~~~\left.-u(s)zM(\cos s,\sin s)\sum\limits_{i+j=0}^3b_{i,j}^-u(s)^{i+2j}z^{i+2j}\cos^is\sin^js\right]{\rm d}s.
\end{split}}$$
Denote by $h(z)=z^3h_1(z)$, which has  the same number of zeros with $h_1(z)$  on $(0,+\infty)$. Then, by direct computation, we have that
\begin{align*}
h(z)&=z^3h_1(z)\\
&=\int_{-\frac{\pi}{2}}^{\frac{\pi}{2}}\frac{1+\sin^2s}{u(s)^4g(s)^2}\left[N(\cos s,\sin s)\sum\limits_{i+j=0}^3a_{i,j}^+u(s)^{i+2(j+1)}z^{i+2(j+1)}\cos^is\sin^js\right.\\
&~~~\left.-M(\cos s,\sin s)\sum\limits_{i+j=0}^3b_{i,j}^+u(s)^{(i+1)+2j}z^{(i+1)+2j}\cos^is\sin^js\right]{\rm d}s\\
&+\int_{\frac{\pi}{2}}^{\frac{3\pi}{2}}\frac{1+\sin^2s}{u(s)^4g(s)^2}\left[N(\cos s,\sin s)\sum\limits_{i+j=0}^3a_{i,j}^-u(s)^{i+2(j+1)}z^{i+2(j+1)}\cos^is\sin^js\right.\\
&~~~\left.-M(\cos s,\sin s)\sum\limits_{i+j=0}^3b_{i,j}^-u(s)^{(i+1)+2j}z^{(i+1)+2j}\cos^is\sin^js\right]{\rm d}s\\
&=\sum\limits_{i+j=1}^4k_{i,j}z^{i+2j}\tag{3.6}
\end{align*}
where $k_{i,j}$'s are some real constants. Specifically,
$$k_{i,j}=a_{i,j-1}^+\int_{-\frac{\pi}{2}}^{\frac{\pi}{2}}\varphi_{i,j}(s){\rm d}s-b_{i-1,j}^+ \int_{-\frac{\pi}{2}}^{\frac{\pi}{2}}\psi_{i,j}(s){\rm d}s+a_{i,j-1}^-\int_{\frac{\pi}{2}}^{\frac{3\pi}{2}}\varphi_{i,j}(s){\rm d}s
-b_{i-1,j}^- \int_{\frac{\pi}{2}}^{\frac{3\pi}{2}}\psi_{i,j}(s){\rm d}s,\eqno{(3.7)}$$

\begin{align*}
\varphi_{i,j}(s)&=\frac{1}{g(s)^2}(1+\sin^2s)N(\cos s,\sin s)\cos^is\sin^{j-1}su(s)^{i+2j-4},\tag{3.8}\\
\psi_{i,j}(s)&=\frac{1}{g(s)^2}(1+\sin^2s)M(\cos s,\sin s)\cos^{i-1}s\sin^jsu(s)^{i+2j-4},\tag{3.9}
\end{align*}
and $a_{i,-1}^{\pm}=b_{-1,j}^{\pm}=0$. Therefore, by $(3.6)$, we know that
$${\begin{split}
h(z)&=z^3h_1(z)\\
&=k_{1,0}z+(k_{0,1}+k_{2,0})z^2+(k_{1,1}+k_{3,0})z^3+(k_{0,2}+k_{2,1}+k_{4,0})z^4+(k_{1,2}+k_{3,1})z^5\\
&+(k_{0,3}+k_{2,2})z^6+k_{1,3}z^7+k_{0,4}z^8.
\end{split}}$$
Hence, using the Theorem 2.2, the first order averging function $h_1(z)$ has at most $7$ positive simple zeros which provide 7 limit cycles of system $(1.4)$, when the averaged function is non-zero. This ends the proof of statement (a) of Theorem 1.1.

Next, we will prove the statement (b)  of Theorem 1.1 by providing a specific example.

Taking $a=b=1,~c=-\frac{1}{4},~d=3$ in $(1.3)$, we obtain the cubic polynomial differential system with a global center at the origin
$$
\begin{cases}
\displaystyle\frac{{\rm d}x}{{\rm d}t}=x^2+y,\\
\displaystyle\frac{{\rm d}y}{{\rm d}t}=-\frac{x^3}{4}+3xy.
\end{cases}\eqno{(3.10)}
$$
Consider the perturbation of system $(3.10)$ as
$$
\begin{cases}
\displaystyle\frac{{\rm d}x}{{\rm d}t}=x^2+y+\varepsilon p_1(x,y),\\
\displaystyle\frac{{\rm d}y}{{\rm d}t}=-\frac{x^3}{4}+3xy+\varepsilon q_1(x,y).
\end{cases}\eqno{(3.11)}
$$
where
$$p_1(x,y)=\begin{cases}a_{0,0}^++a_{1,0}^+x+a_{2,0}^+x^2+a_{3,0}^+x^3+a_{2,1}^+x^2y+a_{1,2}^+xy^2+a_{0,3}^+y^3,~~x\geq0,\\
0,~~x<0,
\end{cases}$$
$$q_1(x,y)=\begin{cases}b_{0,0}^+,~x\geq0,\\
0,~~x<0.
\end{cases}$$
The weighted blow-up $x=r\cos\theta,~y=r^2\sin\theta$, transforms the perturbed system $(3.11)$ into the standard form $(3.1)$ with
$$F_0(\theta,r)=\frac{\cos^3\theta+\cos\theta\sin\theta-\frac{1}{4}\sin\theta\cos^3\theta+3\cos\theta\sin^2\theta}{-\frac{1}{4}\cos^4\theta+\cos^2\theta\sin\theta-2\sin^2\theta},$$
$${\begin{split}
F_1(\theta,r)&=\frac{1+\sin^2\theta}{r^3(-\frac{1}{4}\cos^4\theta+\cos^2\theta\sin\theta-2\sin^2\theta)^2}\left[r^2(-\frac{1}{4}\cos^3\theta+3\cos\theta\sin\theta)p_1(r\cos\theta,r^2\sin\theta)\right.\\
&-\left.r(\cos^2\theta+\sin\theta)q_1(r\cos\theta,r^2\sin\theta)\right].
\end{split}}$$
And in this case, the solution of equation
$$\frac{{\rm d}r}{{\rm d}\theta}=F_0(\theta,r),$$
with initial vale $r(-\frac{\pi}{2})=z$ is
$$r_0(\theta;z)=u(\theta)z,$$
where
$$u(\theta)=\frac{2^{\frac{3}{2}}{\rm exp}\left(-\frac{5}{2}\arctan\frac{4\sin\theta}{-4\sin\theta+\cos(2\theta)+1}\right)}{{\rm exp}(\frac{5\pi}{8})\sqrt[4]{-8\sin\theta-8\sin(3\theta)-28\cos(2\theta)+\cos(4\theta)+35}}.$$
Thus, the first order averaging function is
$${\begin{split}
h_1(z)&=\int_{-\frac{\pi}{2}}^{\frac{\pi}{2}}\frac{F_1(s,u(s)z)}{u(s)}{\rm d}s\\
&=\frac{1}{z^3}\int_{-\frac{\pi}{2}}^{\frac{\pi}{2}}\frac{1+\sin^2s}{u(s)^4(-\frac{1}{4}\cos^4s+\cos^2s\sin s-2\sin^2s)^2}\left[-b_{0,0}^+A_1(s)z+a_{0,0}^+A_2(s)z^2\right.\\
&~~~\left.+a_{1,0}^+A_3(s)z^3+a_{2,0}^+A_4(s)z^4+a_{3,0}^+A_5(s)z^5+a_{2,1}^+A_6(s)z^6+a_{1,2}^+A_7(s)z^7+a_{0,3}^+A_8(s)z^8\right]{\rm d}s
\end{split}}$$
where
$${\begin{split}
A_1(s)&=\left(\cos^2s+\sin s\right)u(s)\\
A_2(s)&=\left(-\frac{1}{4}\cos^3s+3\cos s\sin s\right)u(s)^2,\\
A_3(s)&=\left(-\frac{1}{4}\cos^4s+3\cos^2 s\sin s\right)u(s)^3,\\
A_4(s)&=\left(-\frac{1}{4}\cos^5s+3\cos^3 s\sin s\right)u(s)^4,\\
A_5(s)&=\left(-\frac{1}{4}\cos^6s+3\cos^4 s\sin s\right)u(s)^5,\\
A_6(s)&=\left(-\frac{1}{4}\cos^5s\sin s+3\cos^3 s\sin^2 s\right)u(s)^6,\\
A_7(s)&=\left(-\frac{1}{4}\cos^4s\sin^2s+3\cos^2 s\sin^3 s\right)u(s)^7,\\
A_8(s)&=\left(-\frac{1}{4}\cos^3s\sin^3s+3\cos s\sin^4 s\right)u(s)^8.
\end{split}}$$
By numerical calculation, we can obtain that
$${\begin{split}
h_1(z)&=\frac{1}{z^3}(-15489718.20..b_{0,0}^+z+82848.95524..a_{0,0}^+z^2+740.4727979..a_{1,0}^+z^3\\
&~~~+12.56637060..a_{2,0}^+z^4+24.91789286..a_{3,0}^+z^5+114.4398363..a_{2,1}^+z^6\\
&~~~+540.9497062..a_{1,2}^+z^7+2670.453320..a_{0,3}^+z^8).
\end{split}}$$
Taking
$${\begin{split}
&b_{0,0}^+=\frac{5040}{15489718.20..},~~a_{0,0}^+=\frac{13068}{82848.95524..},~~a_{1,0}^+=\frac{-13132}{740.4727979..}\\
&a_{2,0}^+=\frac{6769}{12.56637060..},~~a_{3,0}^+=\frac{-1960}{24.91789286..},~~a_{2,1}^+=\frac{322}{114.4398363..},\\
&a_{1,2}^+=\frac{-28}{540.9497062..},~~a_{0,3}^+=\frac{1}{2670.453320..},
\end{split}}$$
the function $h_1(z)$ becomes
$$h_1(z)=\frac{1}{z^3}(-5040z+13068z^2-13132z^3+6769z^4-1960z^5+322z^6-28z^7+z^8).$$
It is easy to check that $h_1(z)$ has exactly 7 positive simple zeros given by $z_i=i,~~i=1,\cdots,7$, which provide 7 limit cycles of the perturbed system $(3.11)$ with $\varepsilon\neq0$ sufficiently small. Thus, the statement (b) of Theorem 1.1 is proved.

\vskip 0.2 true cm
\centerline {\bf \large 4 Proof of Theorem 1.2}
\vskip 0.2 true cm

Taking the weighted blow-up $x=r\cos\theta,~y=r\sin^2\theta$ for system $(1.5)$, we can obtain the standard form
$$\frac{{\rm d}r}{{\rm d}\theta}=F_0(\theta,r)+\varepsilon F_1(\theta,r)+O(\varepsilon),\eqno{(4.1)}$$
where
$$F_0(\theta,r)=\frac{f(\theta)}{g(\theta)}r,\eqno{(4.2)}$$
$${\begin{split}
F_1(\theta,r)&=\frac{1+\sin^2\theta}{r^3g(\theta)^2}\left[r^2N(\cos\theta,\sin\theta)p_2(r\cos\theta,r^2\sin\theta)\right.\\
&-\left.rM(\cos\theta,\sin\theta)q_2(r\cos\theta,r^2\sin\theta)\right],
\end{split}}\eqno{(4.3)}$$
and
$$p_2(r\cos\theta,r^2\sin\theta)=\begin{cases}\sum\limits _{i+j=0}^3c_{i,j}^+r^{i+2j}\cos^i\theta \sin^j\theta,~~\sin\theta\geq0,\\
\sum\limits _{i+j=0}^3c_{i,j}^-r^{i+2j}\cos^i\theta \sin^j\theta,~~\sin\theta<0,
\end{cases}$$
$$q_2(r\cos\theta,r^2\sin\theta)=\begin{cases}\sum\limits _{i+j=0}^3d_{i,j}^+r^{i+2j}\cos^i\theta \sin^j\theta,~~\sin\theta\geq0,\\
\sum\limits _{i+j=0}^3d_{i,j}^-r^{i+2j}\cos^i\theta \sin^j\theta,~~\sin\theta<0,
\end{cases}$$
$M,N,f,g$ are given in  $(2.6)-(2.9)$, respectively.

The solution of equation $(4.1)_{\varepsilon=0}$ with initial value $r(0)=z$ is
$$r_0(\theta;z)=v(\theta)z,\eqno{(4.4)}$$
where $v(\theta)={\rm exp}\left(\int_0^{\theta}\frac{f(\tau)}{g(\tau)}{\rm d}\tau\right)$. From $(4.3)$ we know that $F_1(\theta,r)$ is analytic function $2\pi-$periodic in variable $\theta$. By $(2.5)$, we can obtain the first order averaging function of system $(4.1)$
$${\begin{split}
h_1(z)&=\int_{0}^{2\pi}\frac{F_1(s,v(s)z)}{v(s)}{\rm d}s\\
&=\int_{0}^{\pi}\frac{1+\sin^2\theta}{v(s)^4z^3g(s)^2}\left[v(s)^2z^2N(\cos s,\sin s)\sum\limits_{i+j=0}^3c_{i,j}^+v(s)^{i+2j}z^{i+2j}\cos^is\sin^js\right.\\
&~~~\left.-v(s)zM(\cos s,\sin s)\sum\limits_{i+j=0}^3d_{i,j}^+v(s)^{i+2j}z^{i+2j}\cos^is\sin^js\right]{\rm d}s\\
&+\int_{\pi}^{2\pi}\frac{1+\sin^2\theta}{v(s)^4z^3g(s)^2}\left[v(s)^2z^2N(\cos s,\sin s)\sum\limits_{i+j=0}^3c_{i,j}^-v(s)^{i+2j}z^{i+2j}\cos^is\sin^js\right.\\
&~~~\left.-v(s)zM(\cos s,\sin s)\sum\limits_{i+j=0}^3d_{i,j}^-v(s)^{i+2j}z^{i+2j}\cos^is\sin^js\right]{\rm d}s.
\end{split}}$$
Hence, the number of zeros of  function $h_1(z)$ is equal to that of function $\bar{h}(z)=z^3h_1(z)$  on $(0,+\infty)$. Then, by direct computation, we have that
\begin{align*}
\bar{h}(z)&=z^3h_1(z)\\
&=\int_{0}^{\pi}\frac{1+\sin^2s}{v(s)^4g(s)^2}\left[N(\cos s,\sin s)\sum\limits_{i+j=0}^3c_{i,j}^+v(s)^{i+2(j+1)}z^{i+2(j+1)}\cos^is\sin^js\right.\\
&~~~\left.-M(\cos s,\sin s)\sum\limits_{i+j=0}^3d_{i,j}^+v(s)^{(i+1)+2j}z^{(i+1)+2j}\cos^is\sin^js\right]{\rm d}s\\
&+\int_{\pi}^{2\pi}\frac{1+\sin^2s}{v(s)^4g(s)^2}\left[N(\cos s,\sin s)\sum\limits_{i+j=0}^3c_{i,j}^-v(s)^{i+2(j+1)}z^{i+2(j+1)}\cos^is\sin^js\right.\\
&~~~\left.-M(\cos s,\sin s)\sum\limits_{i+j=0}^3d_{i,j}^-v(s)^{(i+1)+2j}z^{(i+1)+2j}\cos^is\sin^js\right]{\rm d}s\\
&=\sum\limits_{i+j=1}^4l_{i,j}z^{i+2j}\tag{4.5}
\end{align*}
where $l_{i,j}$'s are some real constants. Specifically,
$$l_{i,j}=c_{i,j-1}^+\int_{0}^{\pi}\bar{\varphi}_{i,j}(s){\rm d}s-d_{i-1,j}^+ \int_{0}^{\pi}\bar{\psi}_{i,j}(s){\rm d}s+c_{i,j-1}^-\int_{\pi}^{2\pi}\bar{\varphi}_{i,j}(s){\rm d}s
-d_{i-1,j}^- \int_{\pi}^{2\pi}\bar{\psi}_{i,j}(s){\rm d}s,\eqno{(4.6)}$$

\begin{align*}
\bar{\varphi}_{i,j}(s)&=\frac{1}{g(s)^2}(1+\sin^2s)N(\cos s,\sin s)\cos^is\sin^{j-1}sv(s)^{i+2j-4},\tag{4.7}\\
\bar{\psi}_{i,j}(s)&=\frac{1}{g(s)^2}(1+\sin^2s)M(\cos s,\sin s)\cos^{i-1}s\sin^jsv(s)^{i+2j-4},\tag{4.8}
\end{align*}
and $c_{i,-1}^{\pm}=d_{-1,j}^{\pm}=0$.

\noindent{\bf Lemma 4.1.} If $i$ is even, it holds that
$$\bar{\varphi}_{i,j}(\pi-s)=-\bar{\varphi}_{i,j}(s),~~~\bar{\psi}_{i,j}(\pi-s)=-\bar{\psi}_{i,j}(s).\eqno{(4.9)}$$

\noindent{\bf Proof.} We only prove the first equation of $(4.9)$, one can obtain the second one in a similar way.

Note that $\cos(\pi-s)=-\cos s,~\sin(\pi-s)=\sin s$. Form $(2.6)-(2.9)$, we have that
$${\begin{split}
f(\pi-s)&=-\cos sM(-\cos s,\sin s)+\sin sN(-\cos s,\sin s)\\
&=-\cos sM(\cos s,\sin s)-\sin sN(\cos s,\sin s)=-f(s)\end{split}}\eqno{(4.10)}$$
and
 $$g(\pi-s)=g(s)\eqno{(4.11)}$$
analogously. By $(4.10)$ and $(4.11)$, taking the change $s=\pi-\tau$, we get that
$$\int_{0}^{\pi-\theta}\frac{f(s)}{g(s)}{\rm d}s=\int_{\pi}^{\theta}\frac{f(\pi-\tau)}{g(\pi-\tau)}(-{\rm d}\tau)=\int_{\pi}^{\theta}\frac{f(\tau)}{g(\tau)}{\rm d}\tau,\eqno{(4.12)}$$
And doing $s=\pi-\tau$, we have that
$$\int_{0}^{\pi}\frac{f(s)}{g(s)}{\rm d}s=\int_{\pi}^0\frac{f(\pi-\tau)}{g(\pi-\tau)}({-\rm d}\tau)=-\int_0^\pi\frac{f(\tau)}{g(\tau)}{\rm d}\tau,$$
which imply that
$$\int_0^{\pi}\frac{f(\tau)}{g(\tau)}{\rm d}\tau=0.\eqno{(4.13)}$$
Hence, using $(4.12)$ and $(4.13)$, we have that
$$\int_0^{\pi-\theta}\frac{f(s)}{g(s)}{\rm d}s=\int_{\pi}^{0}\frac{f(s)}{g(s)}{\rm d}s+\int_{0}^{\theta}\frac{f(s)}{g(s)}{\rm d}s=\int_{0}^{\theta}\frac{f(s)}{g(s)}{\rm d}s.\eqno{(4.14)}$$
So, by $(4.14)$, we know that
$$v(\pi-\theta)=v(\theta).\eqno{(4.15)}$$
Using $(4.11)$ and $(4.15)$, and noting that $i$ is even, we have that
$${\begin{split}
\bar{\varphi}_{i,j}(\pi-s)&=\frac{1+\sin^2s}{g(\pi-s)^2}(-N(\cos s,\sin s))(-1)^i\cos^is\sin^{j-1}sv(\pi-s)^{i+2j-4}\\
&=-\frac{1+\sin^2s}{g(s)^2}N(\cos s,\sin s)\cos^is\sin^{j-1}sv(s)^{i+2j-4}\\
&=-\bar{\varphi}_{i,j}(s),
\end{split}}$$
which ends the proof.

\noindent{\bf Lemma 4.2.} $\bar{\varphi}_{i,j}(s)$ and $\bar{\psi}_{i,j}(s)$ are both  $2\pi-$periodic functions for all $i,j$.

\noindent{\bf Proof.} By $(2.9)$, $(4.7)$ and $(4.8)$, we only need to prove that
$$v(s+2\pi)=v(s).\eqno{(4.16)}$$
From Lemma 1 of  \cite{GLV17}, we know that
$$\int_{0}^{2\pi}\frac{f(\tau)}{g(\tau)}{\rm d}\tau=0.\eqno{(4.17)}$$
Hence, by $(4.17)$, we get
$$\int_{0}^{2\pi+s}\frac{f(\tau)}{g(\tau)}{\rm d}\tau=\int_{0}^{2\pi}\frac{f(\tau)}{g(\tau)}{\rm d}\tau+\int_{2\pi}^{2\pi+s}\frac{f(\tau)}{g(\tau)}{\rm d}\tau=\int_{0}^{s}\frac{f(t+2\pi)}{g(t+2\pi)}{\rm d}t=\int_{0}^{s}\frac{f(\tau)}{g(\tau)}{\rm d}\tau.\eqno{(4.18)}$$
So, we can obtain $(4.16)$ by $(4.4)$ and $(4.18)$.

\noindent{\bf Lemma 4.3.} If $i$ is even, it holds that
$$\int_{0}^{\pi}\bar{\varphi}_{i,j}(s){\rm d}s=\int_{\pi}^{2\pi}\bar{\varphi}_{i,j}(s){\rm d}s=0,\eqno{(4.19)}$$
$$\int_{0}^{\pi}\bar{\psi}_{i,j}(s){\rm d}s=\int_{\pi}^{2\pi}\bar{\psi}_{i,j}(s){\rm d}s=0.\eqno{(4.20)}$$

\noindent{\bf Proof.} Doing the change $s=\pi-t$, and using Lemma 4.1, we have that
$$\int_{0}^{\pi}\bar{\varphi}_{i,j}(s){\rm d}s=\int_{\pi}^{0}\bar{\varphi}_{i,j}(\pi-t)(-{\rm d}t)=\int_{0}^{\pi}\bar{\varphi}_{i,j}(\pi-t){\rm d}t=-\int_{0}^{\pi}\bar{\varphi}_{i,j}(t){\rm d}t,$$
which imply that
$$\int_{0}^{\pi}\bar{\varphi}_{i,j}(s){\rm d}s=0.\eqno{(4.21)}$$
Moreover taking the change $s=\pi-t$, and by Lemma 4.1, we have that
$$\int_{\pi}^{2\pi}\bar{\varphi}_{i,j}(s){\rm d}s=\int_{0}^{-\pi}\bar{\varphi}_{i,j}(\pi-t)(-{\rm d}t)=-\int_{-\pi}^{0}\bar{\varphi}_{i,j}(t){\rm d}t.$$
Let $\tau=t+2\pi$, by Lemma 4.2, we get that
$$\int_{-\pi}^{0}\bar{\varphi}_{i,j}(t){\rm d}t=\int_{\pi}^{2\pi}\bar{\varphi}_{i,j}(\tau-2\pi){\rm d}\tau=\int_{\pi}^{2\pi}\bar{\varphi}_{i,j}(\tau){\rm d}\tau.$$
Combining the above two expressions, we can obtain that
$$\int_{\pi}^{2\pi}\bar{\varphi}_{i,j}(s){\rm d}s=0.\eqno{(4.22)}$$
We have proved $(4.19)$ by $(4.21)$ and $(4.22)$. One can get $(4.20)$  in a similar way.

 Therefore, by $(4.5)$ and Lemma 4.3, we know that
$${\begin{split}
\bar{h}(z)&=z^3h_1(z)\\
&=l_{1,0}z+(l_{1,1}+l_{3,0})z^3+(l_{1,2}+l_{3,1})z^5+l_{1,3}z^7.
\end{split}}$$
Hence, using the Theorem 2.2, the first order averaging function $h_1(z)$ has at most $3$ positive simple zeros which provide 3 limit cycles of system $(1.4)$, when the first order averaging function is non-zero. This ends the proof of statement (a) of Theorem 1.2.

Next, we will prove the  conclusion  $(b)$ of Theorem 1.2.

Consider  system $(3.10)$ with the following perturbations
$$
\begin{cases}
\displaystyle\frac{{\rm d}x}{{\rm d}t}=x^2+y+\varepsilon p_2(x,y),\\
\displaystyle\frac{{\rm d}y}{{\rm d}t}=-\frac{x^3}{4}+3xy+\varepsilon q_2(x,y).
\end{cases}\eqno{(4.23)}
$$
where
$$p_2(x,y)=\begin{cases}c_{1,0}^+x+c_{1,1}^+xy+c_{1,2}^+y^2,~~y\geq0,\\
0,~~y<0,
\end{cases}$$
$$q_2(x,y)=\begin{cases}d_{0,0}^+,~y\geq0,\\
0,~~y<0.
\end{cases}$$

Doing the change $x=r\cos\theta,~y=r^2\sin\theta$, we can  transforms the perturbed system $(4.23)$ into the standard form $(4.1)$ with
$$F_0(\theta,r)=\frac{\cos^3\theta+\cos\theta\sin\theta-\frac{1}{4}\sin\theta\cos^3\theta+3\cos\theta\sin^2\theta}{-\frac{1}{4}\cos^4\theta+\cos^2\theta\sin\theta-2\sin^2\theta},$$
$${\begin{split}
F_1(\theta,r)&=\frac{1+\sin^2\theta}{r^3(-\frac{1}{4}\cos^4\theta+\cos^2\theta\sin\theta-2\sin^2\theta)^2}\left[r^2(-\frac{1}{4}\cos^3\theta+3\cos\theta\sin\theta)p_2(r\cos\theta,r^2\sin\theta)\right.\\
&-\left.r(\cos^2\theta+\sin\theta)q_2(r\cos\theta,r^2\sin\theta)\right],
\end{split}}.$$
And in this case, the solution of equation
$$\frac{{\rm d}r}{{\rm d}\theta}=F_0(\theta,r),$$
with initial vale $r(0)=z$ is
$$r_0(\theta;z)=v(\theta)z,$$
where
$$v(\theta)=\frac{2^{\frac{3}{4}}{\rm exp}\left(-\frac{5}{2}\arctan\frac{4\sin\theta}{-4\sin\theta+\cos(2\theta)+1}\right)}{\sqrt[4]{-8\sin\theta-8\sin(3\theta)-28\cos(2\theta)+\cos(4\theta)+35}}.$$
Therefore, the first order averaging function is
$${\begin{split}
h_1(z)&=\int_{0}^{\pi}\frac{F_1(s,v(s)z)}{v(s)}{\rm d}s\\
&=\frac{1}{z^3}\int_{0}^{\pi}\frac{1+\sin^2s}{v(s)^4(-\frac{1}{4}\cos^4s+\cos^2s\sin s-2\sin^2s)^2}\left[-d_{0,0}^+B_1(s)z+c_{1,0}^+B_2(s)z^3\right.\\
&~~~\left.+c_{1,1}^+B_3(s)z^5+c_{1,2}^+B_4(s)z^7\right]{\rm d}s
\end{split}}$$
where
$${\begin{split}
B_1(s)&=\left(\cos^2s+\sin s\right)v(s)\\
B_2(s)&=\left(-\frac{1}{4}\cos^4s+3\cos^2 s\sin s\right)v(s)^3,\\
B_3(s)&=\left(-\frac{1}{4}\cos^4s\sin s+3\cos^2 s\sin^2 s\right)v(s)^5,\\
B_4(s)&=\left(-\frac{1}{4}\cos^4s\sin^2s+3\cos^2 s\sin^3 s\right)v(s)^7.
\end{split}}$$
By numerical computation, we can obtain that
$${\begin{split}
h_1(z)&=\frac{1}{z^3}(-407552.3744..d_{0,0}^+z+351.8642184..c_{1,0}^+z^3+138.4955380..c_{1,1}^+z^5\\
&~~~+82260.86314..c_{1,2}^+z^7).
\end{split}}$$
Taking
$${\begin{split}
&d_{0,0}^+=\frac{6}{407552.3744..},~~~~~c_{1,0}^+=\frac{11}{351.8642184..},\\
&c_{1,1}^+=\frac{-6}{138.4955380..},~~~~~c_{1,2}^+=\frac{1}{82260.86314..},
\end{split}}$$
we have that
$$h_1(z)=\frac{1}{z^3}(-6z+11z^3-6z^5+z^7).$$
It is easy to check that $h_1(z)$ has exactly 3 positive simple zeros given by $z_i=\sqrt{i},~~i=1,\cdots,3$, which provide 3 limit cycles of the perturbed system $(4.23)$ with $\varepsilon\neq0$ sufficiently small. Hence, we finish  the proof of statement (b) of Theorem 1.2.

\vskip 0.4 true cm

\noindent{\bf Acknowledgements} This work is supported by National Natural Science Foundation of China (Grant Nos. 12101540).
\vskip 0.2 true cm


\begin{thebibliography}{30}

\bibitem{BZ64}{Berezin, I. S.,   Zhidkov, N. P.: Computing Methods, II. Pergamon Press, Oxford (1964)}

\bibitem{L06}{Chen, F., Li,  C.,  Llibre,  J.,  Zhang, Z.:  A unified proof on the weak Hilbert 16th problem for n=2. J. Differential Equations. {\bf 221}, 309--342 (2006) }

\bibitem{dBCK08}{ di Bernardo, M.,  Budd, C. J.,  Champneys,  A. R., Kowalczyk,  P.: Piecewise-Smooth Dynamical Systems: Theory and Applications. Springer-Verlag, London (2008)}

\bibitem{D00} { Hilbert, D.: Mathematische problem. Bul. Amer. Math. Soc. {\bf 8}, 437--479 (1902) }

\bibitem{GLP13}{ Garc\'{i}a, B., Llibre, J.,   P\'{e}rez del R\'{i}o,  J.: Planar quasi-homogeneous polynomial differential systems and their integrability. J. Differential Equations. {\bf 255}, 3185--3204  (2013) }

\bibitem{G03}{ Garc\'{i}a, I.:  On the integrability of quasihomogeneous and related planar vector fields. Internat. J. Bifur. Chaos Appl. Sci. Engrg. {\bf 13},  995--1002  (2003)}

\bibitem{GGL15}{  Gin\'{e}, J.,  Grau, M.,  Llibre, J.: Limit cycles bifurcating from planar polynomial quasi-homogeneous centers. J. Differential Equations. {\bf 259}, 7135--7160 (2015)}

\bibitem{GLV17}{Gin\'{e}, J.,   Llibre, J.,  Valls, C.:  Centers of weight-homogeneous polynomial vector fields on the plane. Proc. Amer. Math. Soc. {\bf 145},  2539--2555 (2017)}

\bibitem{LL16} { Llibre, J.,  Lopes, B. D., de Moraes, J. R.:  Limit cycles bifurcating from the periodic annulus of the weight-homogeneous polynomial centers of weight-degree 2. Appl. Math. Comput. {\bf 274}, 47--54  (2016) }

\bibitem{LM14}{Llibre, J.,  Mereu,  A.:  Limit cycles for discontinuous quadratic differential systems with two zones. J. Math. Anal. Appl. {\bf 413},  763--775 (2014) }

\bibitem{LP09}{ Llibre, J.,  Pessoa, C.:  On the centers of the weight-homogeneous polynomial vector fields on the plane. J. Math. Anal. Appl. {\bf 359}, 722--730 (2009) }

\bibitem{L03} { Li, J.:  Hilbert's 16th problem and bifurcations of planar polynomial vector fields. Internat. J. Bifur. Chaos Appl. Sci. Engrg. {\bf 13}, 47--106 (2003)  }

\bibitem{LCZ17}{  Li, S., Cen,  X.,  Zhao,  Y.: Bifurcation of limit cycles by perturbing piecewise smooth integrable non-Hamiltonian systems. Nonlinear Anal. Real World Appl. {\bf 34},  140--148 (2017)}

\bibitem{LLYZ09}{ Li, W., Llibre, J.,   Yang, J.,  Zhang, Z.:  Limit cycles bifurcating from the period annulus of quasi-homogeneous centers. J. Dynam. Differential Equations.  {\bf 21}, 133--152  (2009).}

\bibitem{LH10}{Liu, X.,   Han, M.:  Bifurcation of limit cycles by perturbing piecewise Hamiltonian systems. Internat. J. Bifur. Chaos Appl. Sci. Engrg. {\bf 20}, 1379--1390   (2010)
    }

\bibitem{S98} { Smale, S.: Mathematical problems for the next century. Math. Intell. {\bf 20}, 7--15 (1998) }

\bibitem{WZ18}{ Wei, L.,  Zhang, X.:  Averaging theory of arbitrary order for piecewise smooth differential systems and its application. J. Dynam. Differential Equations. {\bf 30},  55--79 (2018)}

\bibitem{YZ18}{  Yang, J.,  Zhao, L.: Bounding the number of limit cycles of discontinuous differential systems by using Picard-Fuchs equations. J. Differential Equations. {\bf 264}, 5734--5757 (2018) }

\end{thebibliography}
\end{document}